\documentclass[english,10pt]{article}

\usepackage{indentfirst}
\usepackage{amsfonts}
\usepackage{amsmath}
\usepackage{latexsym}
\usepackage{delarray}
\usepackage[dvips]{graphics}
\usepackage{amsthm, amssymb}
\usepackage[latin1]{inputenc}

\title{Width of the analyticity strip in space variable of viscous Burgers shockwaves}
\author{C\'edric Lejard \footnote{Universit\'e Paris 6 - Pierre et Marie Curie, Institut de Math\'ematiques de Jussieu, 175 rue du Chevaleret, F75013 PARIS, FRANCE,\texttt{lejard@math.jussieu.fr}}}
\newtheorem*{thm}{Theorem 1}

\newtheorem*{prop}{Proposition 2}
\newtheorem*{prop0}{Proposition 1}
\newtheorem*{lemme}{Lemma}
\newtheorem*{defi}{Definition}
\newtheorem*{prop2}{Proposition 3}
\newtheorem*{def1}{Definition}
\newtheorem*{prop3}{Proposition 4}
\newtheorem*{prop4}{Proposition 5}
\newtheorem*{prop5}{Proposition 6}
\newtheorem*{prop6}{Proposition 7}

\begin{document}

\maketitle

\begin{abstract}
Analytic continuation of viscous shock solution for the generalized Burgers equation with polynomial nonlinear source term is investigated. We show that a pertubated wave recovers its analyticity in the space variable in the strip limited by the first pair of singularities of the wave. 
\end{abstract}

\section{Introduction}

Nonlinear parabolic evolution partial differential equations arise in various fields such as fluids mechanics, traffic flow modelling, mathematical biology and economics. Burgers' equation is one of the most popular of such type of equation, it was first introduced by Bateman [\ref{Ba}] and deeply studied by Burgers [\ref{B}] as a one-dimensional model for turbulence. It received since a constant mathematical interest from which Hopf's study [\ref{H}] may particularly be  outlined. 
Consider the Cauchy problem for generalized Burgers' equation : 
\begin{equation}  \label{GBE}
\begin{split} 
&\frac{\partial f(x,t)}{\partial t} +  \frac{ \partial \Phi(f(x,t))}{\partial x} = \nu \frac{\partial^2 f(x,t)}{\partial x^2} \\ 
&f(x,t=0)=f_0(x) \rightarrow \alpha_{\pm}, x\rightarrow \alpha_{ \pm}, \alpha_+> \alpha_- 
\end{split}
\end{equation} 
If $\Phi$ is a function which satisfies the Gelfand-Oleinik (cf. [\ref{Gel}] and [\ref{Ole}]) entropy condition :
\begin{equation}
\forall u\in [\alpha_-,\alpha_+], \quad c= \frac{\Phi(\alpha_+)-\Phi(\alpha_-)}{\alpha_+-\alpha_-} < \frac{\Phi(u)-\Phi(\alpha_-)}{\alpha_+-\alpha_-}
\end{equation}
this equation enables propagation of shockwaves $\tilde{f}(x-ct+d)$, solving the ordinary differential equation :
\begin{equation}
 -c\tilde{f}(x)+ \Phi (\tilde{f}(x)) +k = \nu \tilde{f}'(x), \quad \tilde{f}(\pm \infty)  = \alpha_{\pm}
\end{equation}
Example of such $\Phi$ are given by concave function on the interval $[\alpha_-,\alpha_+]$.
If in addition, $\Phi$ satisfies the Lax [\ref{Lax}] condition,
\begin{equation}
\Phi'(\alpha_+)<c<\Phi'(\alpha_-)
\end{equation}
the wave is said to be non characteristic, and in this case, the large time asymptotics of the Cauchy problem (\ref{GBE}) are the shifted waves $\tilde{f}(x-ct+d)$. In our study, $\Phi$ is assumed to be polynomial. Its degree will be denoted $n$.
In the particular case $\Phi(u) = -\frac{u^2}{2}$,which will be refered further to as the classical case, the wave has the explicit expression :
\begin{equation}
\tilde{f}(x) = c + \frac{\alpha_+- \alpha_-}{2} \tanh{ \left( (x+d)\frac{\alpha_+-\alpha_-}{4\nu} \right) } 
\end{equation} 
the celerity $c$ and the constant $k$ being constrained by the limits $\tilde{f}(\pm \infty) = \alpha_{\pm}$ : 
\begin{equation}
\begin{split}
& c = \frac{\Phi(\alpha_+)-\Phi (\alpha_-)}{\alpha_+-\alpha_-} \\
& k = \frac{ \alpha_- \Phi( \alpha_+) - \alpha_+ \Phi(\alpha_-) }{\alpha_+ - \alpha_-}
\end{split}
\end{equation}
The shift $d$ may be assumed to be zero up to a space translation of the initial data $f_0$.
The aim of this note is to figure out how the width of the analyticity strip of the solution $f(x,t)$ is matching the width $y_0$($=\frac{2\nu \pi}{\alpha_+-\alpha_-}$ for the classical case) of the analyticity strip of its asymptotic regime $\tilde{f}(x-ct)$.
For Navier-Stokes equations, the analyticity radius of the solutions is known to increase provided that the initial data lies in a functional space that forces it to be zero at infinity. More precisely, 
\begin{itemize}
\item Foias and Temam [\ref{FT}] considered space-periodic solutions and proved that the width of the analyticity strip increases like $ \sqrt{t}$ provided the initial data belongs to the Sobolev space $H^1$ and that the $H^1$-norm of the solution remains uniformly bounded in time. They used Gevrey-norm estimates.
\item Gruj\`ic and Kukavica [\ref{GK}] proved the analyticity of local-in-time solutions with initial data lying in a $L^p(\mathbb{R}^n)$-space and also estimated the growth of the width of the analyticity strip as $\sqrt{t}$.
\item  Lemari\'e-Rieusset [\ref{LR}] has also proved a $\sqrt{t}$-growth for solutions with initial data in $S'(\mathbb{R}^3)$  using techniques from harmonic analysis.
\end{itemize}
Classical Burgers equation may be seen as a 1D version of Navier-Stokes equation. We stress that the width of the analyticity strip does not increase indefinitely if  the overfall $\alpha_+-\alpha_-$ of the shockwave does not vanish.  

\bigskip

The author wishes to thank G. M. Henkin, his thesis supervisor, to have shared with him this quite interesting problem.

\section{Statement of the result}

\begin{thm}
Let $f(x,t)$ be a solution of Burgers' equation with initial condition $f_0(x)$ such that $| f_0(x)-\tilde{f}(x) | =_{x\rightarrow \pm \infty} \mathcal{O}(e^{-\alpha |x|})$ for some $\alpha>0$ and such that $\int_{-\infty}^{\infty} (f_0(x)-\tilde{f}(x)) d x =0$. For all $\varepsilon \in ]0,1]$, there exists $T^*>0$ and $M>0$ such that $f(x,t)$ extends holomorphically in variable $z=x+iy$ in the strip $|y| < \min \left( y_0(1-\varepsilon) , M \sqrt{\nu (t-T^*)} \right) $ for $t>T^*$, where $\pm i y_0$ are the closest singularity of the real axis of the wave $\tilde{f}$ 
\end{thm} 

\noindent
\textbf{Remark 1} : This theorem means that after a transient regime, the width of the analyticity strip of a perturbation of a viscous Burgers shockwave is bounded from below by a non decreasing function of time converging monotonically to the width of the analycity strip of the wave, limited by the closest wave's singularity, which is both an infinity and a branching point of order $n-1$. This singularity is therefore a pole only when $n=2$, in the case of the classical Burgers' equation.\\

\noindent
\textbf{Remark 2} : The proof uses Kato's contraction argument for solving Navier-Stokes equation [\ref{K}], some elements of Sattinger's stability theory of waves [\ref{Sa}] and estimates on the convergence rate in time of the real solution to its asymptotic wave [\ref{IO}],[\ref{Xin}].\\

\noindent
\textbf{Remark 3} : At the limit $\nu=0$, the solutions develop a discontinuity called the inviscid shockwave, this phenomenon can be studied as a result of the loss of analyticity, as proposed by Bessis and Fournier in [\ref{BF}]. Poles are staying away from the real axis so far the viscosity does not vanish. \\

\section{Analytic properties of waves}

The celerity $c$ may assumed to be zero, as $\Phi(u)$ may be replaced by $\Phi(u)-cu$ without loss of generality. The viscosity is brought to be equal to the unity after an appropriate scaling of the solution. Traveling waves are solutions of an ordinary differential equation :
\begin{equation} \label{ode0}
\tilde{f}'=\Phi(\tilde{f})-\Phi(\alpha_{\pm})
\end{equation}
The analyticity of the solution will now be investigated.
\begin{prop0}
Travelling waves solutions $\tilde{f}(x)$ of the Burgers equation are analytic functions of $z=x+iy$ except in a set of isolated points where they admit a singularity which is both an infinity and a branching point of order $n-1$. Moreover, $\tilde{f}$ is analytic on a uniform strip $\lbrace x+iy : |y|<y_0 \rbrace$, where $y_0$ is the imaginary part of the closest singularity of $f$ of the real axis.
\end{prop0}

\begin{proof}
Singularities of solutions of such equations were studied for the first time by Briot and Bouquet $[\ref{BB}]$. A more modern approach may be found in $[\ref{Ince}]$. In the neighbourhood of a regular point $(z_0,f_0) \in \bar{\mathbb{C}} \times \mathbb{C}$, the problem :
\begin{equation} \nonumber
\begin{split}
& \tilde{f}'= \Phi(\tilde{f}) - \Phi(\alpha_{\pm}) \\ 
& \tilde{f}(z_0)=f_0
\end{split}
\end{equation}
is known to  admit powers series solutions around $(z_0,f_0)$ with non vanishing analyticity radius. Consider now a point where $\Phi(\tilde{f})= \infty$, ie $\tilde{f}=\infty$. Let $u=\frac{1}{\tilde{f}}$, which yields :
\begin{equation} \nonumber
\frac{du}{dz} = -u^2 (\Phi(u^{-1})-\Phi(\alpha_{\pm}))
\end{equation}
or,
\begin{equation} \label{ode}
\frac{dz}{du} = \frac{1}{-u^2 (\Phi(u^{-1})-\Phi(\alpha_{\pm}))}
\end{equation}
Where,
\begin{equation} \nonumber
\frac{1}{-u^2 (\Phi(u^{-1})-\Phi(\alpha_{\pm}))} = c u^{n-2}(c_0+c_1 u+ \hdots)
\end{equation}
is an analytic function of $u$ in a neighborhood of $u=0$. Hence $(u=0,z_0)$ is a regular point of the above equation, which therefore admits an analytic solution $z(u)$ in a neighbourhood of $u=0$. Differentiating (\ref{ode}) w.r.t $u$ gives :
\begin{equation} \nonumber
\frac{d^{k} z}{du^k} _{|u=0} = 0 \quad , \quad \forall 1 \leq k \leq n-1
\end{equation}
And consequently :
\begin{equation} \nonumber
z(u)-z_0 = u^{n-1}(\gamma_0+\gamma_1 u + \hdots)
\end{equation}
And finally :
\begin{equation}\nonumber
\tilde{f}(z)=\frac{1}{u(z)}=(z-z_0)^{-\frac{1}{n-1}} \left( \delta_0 + \delta_1 (z-z_0) + \hdots \right)
\end{equation}
Which is well the expected type of singularity.

Integrating (\ref{ode0}), we get that :
\begin{equation}
z_0 = \int_0^{\infty} \frac{dF}{\Phi(F)-\Phi(\alpha_{\pm})}
\end{equation}
The value of this integral depends on the contour from $0$ to infinity :
\begin{equation}
z_0 \in z_0^k + 2i\pi \text{Residue}\left( \frac{1}{\Phi(\tilde{F})-\Phi(\alpha_\pm)} \right) \mathbb{Z}
\end{equation}
where $1\leq k \leq n$.
Singularities are regularly spaced on lines parallel to the imaginary axis and there exists at maximum $n= \text{deg}(\Phi)$ of those lines. This achieves the proof : $y_0$, the imaginary part of the closest singularity of the real axis is positive.

\end{proof}

Define :
\begin{equation} 
\begin{split}
 w(iy) & =  \exp \left( -\frac{1}{2} \int_0^{iy} \Phi'(\tilde{f}(\eta)) d\eta \right) 
 = \exp \left( -\frac{1}{2} \int_0^{iy} \frac{\tilde{f}''(\eta)}{\tilde{f}'(\eta)} d\eta \right)\\
 & = \frac{1}{|\tilde{f}'(iy)|^{\frac{1}{2}}} 
  = \frac{1}{|\Phi(\tilde{f}(iy)) - \Phi(\alpha_+)|^{\frac{1}{2}}}
\end{split}
\end{equation}
Consequently the quantities $|\tilde{f}'(iy)|^{1/2} w(iy)$, $|\Phi(\tilde{f})|^{1/2} w(iy)$ and $|\tilde{f}(iy)|^{n/2} w(iy)$ are bounded function of $y$ in the interval $|y| < y_0$.

\section{Analyticity of solutions}

	\subsection{Weighted norm approach}
	
The function $f(x,t) = \tilde{f}(x) + h(x,t) $ solves (\ref{GBE}) if $h(x,t)$ solves :
\begin{equation} 
\frac{\partial h(x,t)}{\partial t} + \frac{\partial{\Phi'(\tilde{f}(x)) h(x,t)}}{\partial x} - \frac{\partial^2 h(x,t)}{\partial x^2} = -\frac{ \partial R(\tilde{f},h)}{\partial x}  
\end{equation} 
with initial data $h_0(x)= f(x)-\tilde{f}(x)$. $R(\tilde{f},h) = \Phi(\tilde{f} +h) - \Phi(\tilde{f}) - \Phi'(\tilde{f}).h$. The potential $H(x,t)= \int_{-\infty}^x h(\eta,t) d\eta$ solves the partial differential equation :
\begin{equation} 
LH(x,t)=\frac{\partial H(x,t)}{\partial t} + \Phi'(\tilde{f}(x)) \frac{\partial{ H(x,t)}}{\partial x} - \frac{\partial^2 H(x,t)}{\partial x^2} = - R(\tilde{f},\frac{\partial H}{\partial \xi})  
\end{equation} 
with initial data $H_0(x)= \int_{-\infty}^x h_0(\eta) d\eta$, yielding the integral equation :
\begin{equation}
H(x,t) = F_0(x,t) - \int_0^t d\tau \int_{-\infty}^{+ \infty} G(x,t-\tau;\xi) R(\tilde{f},\frac{\partial H}{\partial \xi}) d\xi
\end{equation}
where $F_0 = \int_{-\infty}^{\infty} G(x,t;\xi) H_0(\xi) d\xi $ and $G(x,t;\tau)$ is the Green function of the linearized operator $L$. Differentiating wrt $x$, we get the integral equation solved by $h(x,t)$.
\begin{equation}
h(x,t) = \partial_x F_0(x,t) - \int_0^t d\tau \int_{-\infty}^{+ \infty} G(x,t-\tau;\xi) R(\tilde{f}(\xi),h(\xi,\tau)) d\xi
\end{equation}
Following Sattinger [\ref{Sa}], the Green function $G$ of the operator $L$ in the above integral equation may be expressed as :
\begin{prop}
The Green function $G(x,t;\xi)$ of the linearized operator $L$ has the expression :
\begin{equation}
G(x,t;\xi) = \frac{w(\xi)}{w(x)} K(x,t;\xi) 
\end{equation}
where $w$ stands for the weight function 
\begin{equation}
w(x)= \exp \left( -\frac{1}{2} \int_0^x \Phi'(\tilde{f}(\eta)) d \eta \right)
\end{equation}
and $K(x,t;\xi)$ for the Green function of the operator :
\begin{equation}
\mathcal{M}u(x,t) = \frac{\partial u(x,t) }{\partial t} - \frac{\partial^2 u(x,t)}{\partial x^2} - \frac{1}{2} \left( \tilde{f}'(x) \Phi''(\tilde{f}(x)) + \frac{1}{2} \Phi'(\tilde{f})^2 \right) u(x,t)
\end{equation}
\end{prop}

\begin{proof}
Let $H(x,t)=w^{-1}(x) u(x,t)$. A direct calculation
 yields :
 \begin{equation} \nonumber
 \begin{split}
 & \partial_x H(x,t) = w^{-1}(x) \left( \frac{1}{2} \Phi'(\tilde{f}) u(x,t) + \partial_x u(x,t) \right) \\
 & \partial_{xx} H(x,t) = w^{-1}(x) \left( \left( \frac{1}{4} \Phi'(\tilde{f})^2+ \frac{1}{2} \tilde{f}'(x) \Phi''(\tilde{f}) \right) u(x,t) + \Phi'(\tilde{f}) \partial_x u + \partial_{xx}u \right)
 \end{split}
 \end{equation} 
Combining this two relations, we get that $LH=w^{-1} \mathcal{M} w H$. The function $H(x,t)=w^{-1}(x) u(x,t)$ is therefore a solution of $LH(x,t)=g(x,t)$ if and only if $u(x,t)$ solves the linear equation $\mathcal{M} (wH)=wg$. This completes the proof.
 \end{proof}

\begin{defi}
Denote $m(y,t)=\text{sgn}(y) \min (|y|, M\sqrt{t})$ and $ \Delta_{c} = \lbrace x+iy : |y|<c \rbrace$. Introduce the following functional spaces :
\begin{itemize}
\item the weighted Hardy Space 
	\begin{equation}  \nonumber
	\begin{split}
	& \mathcal{H}^2_{w,a,b}(\Delta_c) = \left\lbrace u \in \mathcal{O}(\Delta_c) : u(\mathbb{R})\subset \mathbb{R}, \quad  \|u\|_{\mathcal{H}^2_{w,a,b}(\Delta_c)} = \right.\\
	& \left. \sup_{0 \leq y <c} \frac{|w(iy)| ^{a}}{(\sqrt{c-|y|)}^{b}} \left( \int_{- \infty}^{\infty} | w(\xi+iy) | . |u(\xi+iy,t)|^2 d\xi \right)^{1/2} < \infty \right\rbrace
	\end{split}
	\end{equation}
	 with norm $\|.\|_{\mathcal{H}^2_{w,a,b}(\Delta_c)}$. 
\item the Banach space $\mathcal{B}_{\varepsilon,T}$ of holomorphic functions $f(.,t), 0<t\leq T$ in the strip $\Delta_{m(y_0(1-\varepsilon),t)})$ with norm $ \| . \|_{\varepsilon,T}$ defined as 
\begin{equation} \nonumber
\| h \|_{\varepsilon,T} = \sup_{0<t \leq T} t^{\frac{1}{4(n-1)}}  \| h(.,t) \|_{\mathcal{H}^2_{w,a_n,b_n}(\Delta_{m(y_0(1-\varepsilon),t)})} < \infty
\end{equation}
where,
\begin{equation}
\begin{split}
&a_n = \frac{2}{n}-\frac{5}{2} \\
&b_n = \frac{n-2}{n-1}
\end{split}
\end{equation}
\end{itemize}
\end{defi}

The forthcomming proposition is also inspired by Sattinger with some modifications:
\begin{itemize}
\item Estimates of the analytic continuation of the Green function $K$ in spaces variables are needed.
\item These estimates are realized in a weighted $L^2$-Hardy norm as introduced in the previous definition and not in a first order weighted Sobolev in supremum norm as in Sattinger's work.
\end{itemize}

\begin{prop2} \label{prop2}
The kernel $K(x,t;\xi)$ is analytic in variable $x$ and $\xi$ in the strip $\Delta_{m(y_0,M\sqrt{t})}$ and satisfies the following estimates :
\begin{equation}
\begin{split}
& \| \int_{-\infty}^{\infty} K(x+iy,t,\xi+i\eta) w(\xi+i\eta) h(\xi+i\eta)^2 d\xi \|_{L^2(dx)} \\
& \leq e^{ \frac{(y-\eta)^2}{4t}} \frac{\text{const}. e^{-\omega t}}{t^{\frac{1}{4}}} \| w^{\frac{1}{2}}. h \|^2_{\mathcal{H}^2_{w,0,0}(\Delta_{m(y_0(1-\varepsilon),t)})} \\
& \| \int_{-\infty}^{\infty} \partial_x K(x+iy,t,\xi+i\eta) w(\xi+i\eta) h(\xi+i\eta)^2  d\xi \|_{L^2(dx)} \\
& \leq e^{ \frac{(y-\eta)^2}{4t}} \frac{\text{const}. e^{-\omega t}}{t^{\frac{3}{4}}} \| w^{\frac{1}{2}}. h \|^2_{\mathcal{H}^2_{w,0,0}(\Delta_{m(y_0(1-\varepsilon),t)})}
\end{split}
\end{equation}
\end{prop2}	

%\noindent
%\textbf{Remark} : In a supremum type Hardy space $\mathcal{H}^{\infty}$ equipped with norm of the form $\| f \|_{\mathcal{H}^{\infty}} = \sup_{|y| <| m(y_0,M\sqrt{t})|} sup_{x\in \mathbb{R}} |w(x+iy)|^{1/2} |f(x+iy)|$, the yielded estimates for the kernel and its derivative behave respectively as $t^{-1/2}$ and $t^{-1}$ as $t$ goes to zero, breaking convergence of the iterative argument for the solving of the nonlinear integral equation. 

\begin{proof}
Consider the eigenfunctions equation for the ordinary linear operator $Au=u''+p(x)$ :
\begin{equation}
u''+p(x)u = \lambda u
\end{equation}
Look for a solution exhibiting exponential dichotomy at $ x \rightarrow \pm \infty$ :
\begin{equation}
\begin{split}
& \varphi_+(x,\lambda) = e^{-\gamma_+(\lambda)x} a_+(x,\lambda) , \quad x \rightarrow \infty\\
& \varphi_-(x,\lambda) = e^{\gamma_-(\lambda)x} a_-(x,\lambda) , \quad x \rightarrow -\infty
\end{split}
\end{equation}
Where,
\begin{equation}
\begin{split}
& \gamma_{\pm}(\lambda) = \sqrt{\lambda - p_{\pm}} \\
& p_{\pm} = \lim_{x\rightarrow \pm \infty} p(x) = \frac{1}{4} \Phi'(\alpha_{\pm})^2 \\
& \bar{p} = \max (p_+,p_-)
\end{split}
\end{equation}
Remark that for $|\arg (\lambda - \bar{p})| < \pi$, 
\begin{equation}
\text{Re}(\gamma_{\pm}(\lambda)) >0
\end{equation}
And consequently, the investigated $\varphi_{\pm}$ decrease exponentially as $x$ goes to $\pm \infty$. The functions $a_+$ and $a_-$ are analytic in $x$ and $\lambda$ and satisfy the following asymptotics :
\begin{equation}
a_{\pm}(x+iy,\lambda) = 1 + \mathcal{O}\left(\frac{1}{\sqrt{|\lambda|}}\right), \quad x \rightarrow \pm \infty
\end{equation}
For instance, $a_+$ is a solution as a function of $x$ of the differential equation :
\begin{equation}
a_+''-2 \gamma_+ a_+' + (p-\bar{p}) a_+ = 0
\end{equation}
which can be reformulated as an integral equation :
\begin{equation}
a_+(x,\lambda) = 1+ \int_{x}^{\infty} \left[p\left(\frac{1-e^{2\gamma_+(x-s)}}{2 \gamma_+}\right)-\bar{p}\right] a_+(s,\lambda) ds
\end{equation}
This equation may be readily complexified in variable $z=x+iy$ :
\begin{equation}
a_+(z,\lambda) = 1  + \int_{z}^{\infty} \left[p\left(\frac{1-e^{2\gamma_+(x-s)}}{2 \gamma_+}\right)-\bar{p}\right] a_+(s+iy,\lambda) ds
\end{equation}
On each real line $y=\text{const}$, the function $p-\bar{p}$ is absolutely integrable and $ \int_{x}^{\infty} |p(s+iy)-\bar{p}| ds$ is uniformly bounded wrt $y$ as $x \rightarrow \infty$. The integral equation may therefore be solved, which yields a $a_ +(z,\lambda)= 1+\mathcal{O}\left( \sqrt{\frac{1}{|\lambda|}} \right)$ (and a $a_-$ satisfying the same estimate) .

Denote $F(x,\xi,\lambda)$ the Green function of the operator $\lambda-A$, it satisfies the following estimates :
\begin{equation}
\begin{split} 
& \| F(x+iy,\xi+i\eta,\lambda) \|_{L^2(dx)} \leq \frac{C_1}{|\lambda|^{3/4}} e^{(C_2+|\gamma_+(\lambda)|)(y-\eta)} \\
& \|\partial_x F(x+iy,\xi+i\eta,\lambda) \|_{L^2(dx)} \leq \frac{C_3}{|\lambda|^{1/4}} e^{(C_2+|\gamma_+(\lambda)|)(y-\eta)}
\end{split}
\end{equation}
Classical spectral theory provides an explicit expression of $F(x,\xi,\lambda)$ : 
\begin{equation}
\begin{split} 
F(x,\xi,\lambda)  = & \frac{1}{W(\lambda)} \varphi_-(\xi) \varphi_+(x),\quad \text{if} \, \xi \leq x \\
& \frac{1}{W(\lambda)} \varphi_-(x) \varphi_+(\xi),\quad \text{if} \, \xi \geq x 
\end{split}
\end{equation}
Where $W(\lambda) = -(\gamma_+(\lambda) + \gamma_-(\lambda)) + \mathcal{O}(1) $ is the wronskian of the eigenfunctions equation (cf [\ref{Sa}]). It is independant of $x$ and $\xi$.
Estimates of $\|F(x+iy,\xi+i\eta,\lambda \|_{L^2(d\xi)}$ are now requested.
\begin{equation}
\begin{split}
& \int_{-\infty}^{\infty} |F(x+iy,\xi + i \eta,\lambda)|^2 d\xi = \frac{1}{W(\lambda)} \left( \int_{-\infty}^x  |\varphi_-(\xi+i\eta) \varphi_+(x+iy)|^2 d\xi \right. \\
& \left. + \int_x^{\infty} | \varphi_-(x+iy) \varphi_+(\xi+i\eta) | ^2 d\xi \right)
\end{split}
\end{equation}
We now estimate the different term as in Sattinger work, using the fact that $\|e^{-\sqrt{|\lambda|} x } \|_{L^2(dx)} = \frac{1}{|\lambda|^{1/4}}$, which saves a factor $\frac{1}{|\lambda|^{1/4}}$ compared to $L^1$-norm. For example, consider $x>0$,
\begin{equation}
\begin{split}
\int_{- \infty}^x  |\varphi_-(\xi+i\eta) & \varphi_+(x+iy)|^2 d\xi \leq e^{ -2\text{Re}(\gamma_+(\lambda))x+2\text{Im}(\gamma_+(\lambda))y} \left( \int_{_\infty}^0 |\varphi_+(\xi+i\eta,\lambda)|^2 d\xi \right. \\
& \left. + \int_0^x |\varphi_+(\xi+i\eta,\lambda)|^2 d \xi \right)
\end{split}
\end{equation}
Where,
\begin{equation}
\begin{split}
&\int_0^{\infty} |\varphi_+(\xi+i\eta)|^2 d \xi  \leq e^{ -2\text{Im}(\gamma_-)\eta} \int_0^{\infty} e^{2\text{Re}(\gamma_-) \xi} \left( 1+ \mathcal{O} \left( \frac{1}{\sqrt{|\lambda|}} \right) \right)^2 d \xi \\
& \leq \frac{1}{2 \text{Re}(\gamma_-(\lambda))} \left( 1+ \mathcal{O}(\frac{1}{\gamma_-(\lambda)})\right)^2  \leq \frac{\text{const}}{\sqrt{|\lambda|}}
\end{split}
\end{equation}
After composition with the square-root, we get the expected power in $\lambda$ ($\frac{1}{W(\lambda)}$ is estimated by $\frac{1}{\sqrt{|\lambda|}}$). Others terms are estimated using the same method. To estimate the norm of derivative, it is sufficient to remark that derivation wrt $x$ of $F(x,\xi,\lambda)$ gives rise to an extra $\sqrt{\lambda}$ factor.
Finally, we need to prove that the kernel of the evolution operator $e^{tA}$ satisfies the expected estimates and that its $\mathcal{H}^2$ norm wrt $\xi + i \eta$ is uniformily bounded in $x+iy$. The kernel $ K(x,t,\xi)$ is represented by the following contour integral :
\begin{equation}
K(x+iy,t,\xi+i\eta) = \frac{1}{2i\pi} \int_{\mathcal{C}} e^{t\lambda} F(x+iy,\xi+i\eta,\lambda) d\lambda
\end{equation}
Where $\mathcal{C}$ is any contour avoiding the essential spectrum of $A$. Sattinger used the contour $\lambda= -\omega + \rho e^{\pm i\delta} $, parametrized by $0 \leq \rho < \infty$. The two parameter $0 < \omega < \bar{p}$ and $\pi / 2<\delta < \pi$ are chosen such that the estimates performed on $F(x,t,\xi)$ hold. The same contour will be used. Using the estimates on $F$ previously established, we get :
\begin{equation}
\begin{split}
\| K(x+iy,t,\xi+i\eta) \|_{L^2(d\xi)} \leq \text{const} \int_{\mathcal{C}} \frac{|e^{\lambda t}|}{|\lambda|^{3/4}}e^{C \sqrt{|\lambda|}(y-\eta)} |d\lambda| 
\end{split}
\end{equation}
In order to estimate the right handside, rewrite it as an integral on the parameter $\rho = |\lambda+\omega|$.
\begin{equation}
\begin{split}
\int_{\mathcal{C}} \frac{|e^{\lambda t}|}{|\lambda|^{3/4}}e^{C \sqrt{|\lambda|}(y-\eta)} |d\lambda|  & \leq \text{const} \int_0^{\infty} \frac{e^{-\omega t} e^{\rho t \cos \delta + \sqrt{\rho} \sin  (\pm \delta)(y-\eta) }}{(\rho \cos \delta - \omega)^2 + \rho^2 \sin^2 \delta)^{3/8}} d\rho \\
& \leq e^{ -\omega t} \int_0^{\infty} \frac{e^{\rho t \cos \delta \pm \sqrt{\rho} \sin |\delta| (y-\eta)}}{\omega^{3/4}+\rho^{3/4} |\sin \delta|^{3/4}} d\rho
\end{split}
\end{equation}
With $\frac{\pi}{2} \leq \delta \leq \pi$, $\cos \delta <0$ and the convergence of the integral is assured. Remark that :
\begin{equation}
\begin{split}
\rho t \cos \delta \pm \sqrt{\rho} \sin |\delta| (y-\eta) & = \cos(\delta) \left( \sqrt{\rho t} \pm \frac{\sin{|\delta|}}{\cos \delta} \frac{(y-\eta)}{2 \sqrt{t}} \right)^2 - \tan|\delta| \cos \delta \frac{(y-\eta)^2}{4t} \\
& =  -| \cos(\delta)| \left( \sqrt{\rho t} \pm \frac{\sin{|\delta|}}{\cos \delta} \frac{(y-\eta)}{2 \sqrt{t}} \right)^2 + \sin|\delta| \frac{(y-\eta)^2}{4t}
\end{split}
 \end{equation}
 We conclude by using Young inequality since the kernel $K(x+iy,t,\xi+i\eta)$ is square-integrable with respect to $x$ and $\xi$, and the $L^2(d\xi)$-norm is uniformily bounded wrt $x$. We use that $ \int_0^{\infty} \frac{e^{-\rho t}}{\rho^{\gamma}}d\rho = \frac{\text{const}}{t^{1-\gamma}}$. The estimates are obtained with $\gamma=3/4$ for $K$ and $\gamma=1/4$ for its derivative.
\end{proof}

	\subsection{Solving the integral equation}

\begin{prop3}
Let $(X,\|.\|)$ be a Banach Space and $B:X \rightarrow X$ a $C^1$ map such that $B(0)=0$ and $\displaystyle{\sup_{\|h\| \leq \frac{1}{1-\sigma} \alpha} ||| DB(h) ||| <\sigma<1}$. Let $a\in X$ such that $\|a\|< \alpha$. The sequence defined by the following iterative scheme : 
\begin{equation}
\begin{split}
& h_{n+1}=a+B(h_n) \\
& h_0=a 
\end{split}
\end{equation}
remains in a ball of center $0$ and radius $\frac{1}{1-\sigma} \alpha$ and has  a limit $h$ in X solving the equation $h=a+B(h)$.
\end{prop3} 

\begin{proof}
Assume that $\|h_n\| \leq \frac{1}{1-\sigma}\alpha$. The mean value inequality yields :
\begin{equation} \nonumber
\|h_{n+1} \| \leq \alpha + \sigma \|h_n \| \leq \alpha + \frac{\sigma}{1-\sigma} \alpha \leq \frac{1}{1-\sigma} \alpha
\end{equation}
The sequence $(h_n)$ is therefore bounded. Using the mean value inequality again, we get :
\begin{equation} \nonumber
\| h_{n+1}- h_n \| \leq \sigma \| h_n - h_{n-1} \|
\end{equation}
Picard's contraction theorem is enough to conclude.
\end{proof}

\begin{lemme}
The evaluation map $\mathcal{H}^2_{w,a,b}(\Delta_c) \rightarrow \mathbb{C}$ is continuous :
\begin{equation}
|u(z)| \leq  \frac{\text{const}.(c-|y|)^{\frac{b-1}{2}}}{|w(iy)|^{\frac{1}{2}+a}} \| u \|_{\mathcal{H}^2_{w,a,b}(\Delta_c)}
\end{equation}
\end{lemme}

\begin{proof}
Represent $w(z)u(z)^2$ as an integral using Cauchy's formula over the contour 
\begin{equation} \nonumber
\mathcal{C}= \lbrace \xi+ic : \xi \in  [\infty,-\infty] \rbrace \cup \lbrace \xi-ic : \xi \in  [-\infty,\infty] \rbrace 
\end{equation}
that is to say,
\begin{equation} 
w(z)u(z)^2 = \frac{1}{2 i \pi} \int_{\mathcal{C}} \frac{w(\zeta)u(\zeta)^2 d\zeta}{\zeta-z}
\end{equation}
It provides the estimate :
\begin{equation}
|u(z)|^2 \leq \frac{1}{\sqrt{\pi}} \frac{1}{c-|y|} \frac{1}{|w(iy)|} \int_{-\infty}^{\infty} |w(\xi+ic)|.|u(\xi+ic)|^2 d \xi
\end{equation}
which is the expected estimate.
\end{proof}

\begin{prop4}
The nonlinear operator $B:\mathcal{B}_{\varepsilon,T} \rightarrow \mathcal{B}_{\varepsilon, T}$,
\begin{equation}
\begin{split}
& B(h)(z,t)= \int_{0}^t d\tau \int_{-\infty}^{\infty} \frac{\partial G}{\partial z}\left(x+iy,t-\tau,\xi+i m(y,\tau) \right) \\
& R \left( \tilde{f}(\xi+im(y,\tau)),h(\xi+im(y,\tau)) \right) d\xi
\end{split}
\end{equation}
with $R(\tilde{f},h) =  \Phi(\tilde{f}+h)-\Phi(\tilde{f}) - \Phi'(\tilde{f}).h $,
satisfies the estimate :
\begin{equation}
\sup_{\|h\|_{\varepsilon,T}\leq 2 \alpha} |||DB(h)||| \leq  C \alpha p(\alpha)
\end{equation}
With $C$ a constant and $p$ a polynomial of degree $\leq n-2$, both being independant of $\varepsilon$ and $T$.
\end{prop4}

\begin{proof}
Differentiating with respect to $h$ yields :
\begin{equation}
\begin{split}
&DB(h)(g)(z,t)= \int_{0}^t d\tau \int_{-\infty}^{\infty} \frac{\partial G}{\partial z}\left(x+iy,\xi+i m(y,\tau)),t-\tau \right) \\
& \frac{\partial{R}}{\partial h} \left( \tilde{f}(\xi+im(y,\tau)),h(\xi+im(y,\tau),\tau) \right) g(\xi+im(y,\tau),\tau) d\xi
\end{split}
\end{equation}
where ,
\begin{equation} \nonumber
\frac{\partial{R}}{\partial h}(\tilde{f},h) = h \sum_{k=0}^{n-2} \frac{1}{(k+1)!} \Phi^{(k+2)}(\tilde{f}) h^{k} 
\end{equation}
For $0 \leq k \leq n-2$, define the following operators :
\begin{eqnarray}
A^{k}_1(h,g)(z,t)& = -& \int_{0}^t \frac{d \tau}{w(z)} \int_{-\infty}^{\infty} w(\xi+im(y,\tau)) \partial_x K \left( (x+iy),(\xi+i m(y,\tau)),t-\tau \right) \nonumber\\
 					 &   &  \Phi^{(k+2)}(\tilde{f}(z))  h^{k+1}(\xi+im(y,\tau),\tau) g(\xi+im(y,\tau),\tau) d\xi \nonumber \\
A^{k}_2(h,g)(z,t) & = & -\frac{1}{2} \int_{0}^t \frac{d \tau}{w(z) } \int_{-\infty}^{\infty} w(\xi+i m(y,\tau))  K \left( x+iy,\xi+i m(y,\tau),t-\tau \right)\nonumber \\
 						&    &   \Phi^{(k+2)}(\tilde{f}(z)) \Phi'(\tilde{f}(z))  h^{k+1}(\xi+im(y,\tau),\tau) g(\xi+im(y,\tau),\tau) d\xi  \nonumber
\end{eqnarray}
%---------------------------------------------------------------------------------------------------------------------------------------------------------------------------------
%---------------------------------------------------------------------------------------------------------------------------------------------------------------------------------
The expected estimates will be proved separately for each $A^{k}_1$ and each $A^{k}_2$. Proceed first with the $A^{k}_1$'s :
\begin{equation} \nonumber
\begin{split}
& \left\| \int_{-\infty}^{\infty} w(\xi+im(y,\tau))  \partial_x K( x+iy,\xi+i m(y,\tau),t-\tau) \Phi^{(k+2)}(\tilde{f}(z)) h^{k} h g(\xi+im(y,\tau),\tau) d\xi \right\|_{L^2(d x)}  \\
&  \leq e^{M^2/4} \text{const} |\Phi^{k+2}(\tilde{f}(z))| . |w(iy)|^{-(\frac{1}{2}+a_n)k}(t-\tau)^{-3/4} (y_0(1-\varepsilon)-|y|)^{k\frac{b_n-1}{2}} \|h  \|^k_{\mathcal{H}^2_w(\Delta_{m(y_0(1-\varepsilon),\tau)})}\\
& \times  \int_{- \infty}^{\infty}  |w(\xi+im(y,\tau)) | . |h(\xi+im(y,\tau),\tau)|. |g(\xi+im(y,\tau),\tau)| d\xi =I_1
\end{split}
\end{equation}
Where the lemma about the continuity of the evaluation map in $\mathcal{H}^2_w(\Delta_{m(y_0(1-\varepsilon),\tau)})$ and the estimate (proposition 2) concerning  the kernel $\partial_z K$ have been used. Moreover, the inequality $\Phi^{(k+2)}(\tilde{f}(z)) \leq \text{const} w(iy)^{\frac{-2}{n}(n-k-2)}$ hold. The Cauchy-Schwarz inequality leads to :
\begin{equation} \nonumber
\begin{split}
I_1  \leq &  e^{M^2/4} \text{const} |w(iy)|^{-\frac{2}{n}(n-k-2)-(\frac{1}{2}+a_n)k-2a_n}(t-\tau)^{-3/4} (y_0(1-\varepsilon)-|y|)^{k\frac{b_n-1}{2}+b_n}  \\
& \times \|h  \|^{k+1}_{\mathcal{H}^2_w(\Delta_{m(y_0(1-\varepsilon),\tau)})} \|g  \|_{\mathcal{H}^2_w(\Delta_{m(y_0(1-\varepsilon),\tau)})}
\end{split}
\end{equation}
Hence,
\begin{equation} \nonumber
\begin{split}
& |||A^{k}_1(h,.) |||_{\mathcal{L}(\mathcal{H}^2)} \leq   e^{M^2/4} \text{const} |w(iy)|^{-\frac{2}{n}(n-k-2)-(\frac{1}{2}+a_n)k-a_n-\frac{1}{2}} \\
 & \times (y_0-|y|)^{k\frac{b_n-1}{2}+\frac{b_n}{2}} \|h  \|^{k+1}_{\varepsilon,t}  \int_{0}^t \frac{e^{-\omega t} d \tau}{(t-\tau)^{\frac{3}{4}} \tau^{\frac{k+2}{4(n-1)}}}  
\end{split}
\end{equation}
Where,
\begin{equation} \nonumber
\int_{0}^t \frac{d \tau}{(t-\tau)^{\frac{3}{4}} \tau^{\frac{k+2}{4(n-1)}}}  = \text{const}t^{1-(\frac{3}{4}+\frac{k+2}{4(n-1)})} =  \text{const} \,  t^{\frac{n-k-3}{4(n-1)}} 
\end{equation}
And finally,
\begin{equation} \nonumber
\begin{split}
\sup_{\|h\|_{\varepsilon,t}\leq 2 \alpha }|||A^{k}_1(h,.) |||_{\varepsilon,t} & \leq  \text{const} e^{\frac{M^2}{4}} |w(iy)|^{-\frac{2}{n}(n-k-2)-(\frac{1}{2}+a_n)k-a_n-\frac{1}{2})} \\ & \times (y_0-|y|)^{k\frac{b_n-1}{2}+\frac{b_n}{2}} e^{-\omega t} t^{\frac{n-k-2}{4(n-1)}} \alpha^{k+1} \\
& \leq \text{const} \alpha^{k+1} e^{-(\omega-\tilde{\omega})t}
\end{split}
\end{equation} 
The last inequality is obtained through the following considerations :
\begin{itemize}
\item The function $y \mapsto |w(iy)|$ is bounded on $[0,y_0]$ and $-\frac{2}{n}(n-k-2)-(\frac{1}{2}+a_n)k-a_n-\frac{1}{2} \geq 0$, with $a_n = \frac{2}{n}-\frac{5}{2}$.
\item The function $y\mapsto y_0-|y|$ is bounded and ${k\frac{b_n-1}{2}+\frac{b_n}{2}} \geq 0$, with $b_n=\frac{n-2}{n-1}$.
\item For $q \geq 0$ and $0\leq \tilde{\omega} \leq \omega$, $t\mapsto t^{q}e^{- \tilde{\omega} t}$ is bounded and $q=\frac{n-k-2}{4(n-1)} \geq 0$.
\end{itemize}
%----------------------------------------------------------------------------------------------------------------------------------------------------------------------------------
%----------------------------------------------------------------------------------------------------------------------------------------------------------------------------------
Proceed now with $A^{k}_2$
\begin{equation} \nonumber
\begin{split}
& \left\| \int_{-\infty}^{\infty} w(\xi+im(y,\tau)) K( x+iy,\xi+i m(y,\tau),t-\tau) \Phi^{(k+2)}(\tilde{f}(z)) \Phi'(\tilde{f}(z)) h^{k} h g(\xi+im(y,\tau),\tau) d\xi \right\|_{L^2(d x)} \\ 
& \leq  e^{M^2/4} \text{const} |\Phi^{(k+2)}(\tilde{f}(z)). \Phi'(\tilde{f}(z))|. |w(iy)|^{-(\frac{1}{2}+a_n)k}(t-\tau)^{-1/4} (y_0-|y|)^{k\frac{b_n-1}{2}} \|h  \|^k_{\mathcal{H}^2_w(\Delta_{m(y_0(1-\varepsilon),\tau)})} \\
& \times \int_{- \infty}^{\infty}  |w(\xi+im(y,\tau)) | . |h(\xi+im(y,\tau),\tau)|. |g(\xi+im(y,\tau),\tau)| d\xi =I_2
\end{split}
\end{equation}
Using an analoguous argument as previously, we get :
\begin{equation} \nonumber
\begin{split}
I_1 \leq &e^{M^2/4} \text{const} |w(iy)|^{-\frac{2}{n}(2n-k-3)-(\frac{1}{2}+a_n)k-2a_n}(t-\tau)^{-1/4} (y_0-|y|)^{k\frac{b_n-1}{2}+b_n}\\
& \times \|h  \|^{k+1}_{\mathcal{H}^2_w(\Delta_{m(y_0(1-\varepsilon),\tau)})} \|g  \|_{\mathcal{H}^2_w(\Delta_{m(y_0(1-\varepsilon),\tau)})}
\end{split}
\end{equation}
And consequently,
\begin{equation} \nonumber
 |||A^{k}_2(h,.) |||_{\mathcal{L}(\mathcal{H}^2)} \leq   e^{M^2/4} \text{const}. |w(iy)|^{-\frac{2}{n}(2n-k-3)-(\frac{1}{2}+a_n)k-a_n-\frac{1}{2}} \|h  \|^{k+1}_{\varepsilon,t}  \int_{0}^t \frac{e^{-t} d \tau}{(t-\tau)^{\frac{1}{4}} \tau^{\frac{k+2}{4(n-1)}}}  
\end{equation}
Where,
\begin{equation} \nonumber
\int_{0}^t \frac{d \tau}{(t-\tau)^{\frac{1}{4}} \tau^{\frac{k+2}{4(n-1)}}}  = \text{const}t^{1-(\frac{1}{4}+\frac{k+2}{4(n-1)})} =  \text{const}. t^{\frac{3n-k-5}{4(n-1)}}  
\end{equation}
It eventually yields the estimate :
\begin{equation} \nonumber
\begin{split}
\sup_{\|h\|_{\varepsilon,t}\leq 2 \alpha }|||A^{k,k'}_2(h,.) |||_{\varepsilon,t} & \leq  \text{const}. e^{\frac{M^2}{4}} |w(iy)|^{-\frac{2}{n}(2n-k-3)-(\frac{1}{2}+a_n)k-a_n-\frac{1}{2}} e^{-\omega t} t^{\frac{3n-k-5}{4(n-1)}} \alpha^{k+1} \\
& \leq \text{const}. \alpha^{k+1} e^{-(\omega-\tilde{\omega})t}
\end{split}
\end{equation} 
The differential $DB(h)(g)$ being a linear superposition of the monomial operators $A^{k}_1(h,g)$ and $A^{k}_2(h,g)$, we get the required estimate.
\end{proof}

\begin{def1}
Introduce the weighted energy $E_H(t)$ of a function $H(.,t)$ :
$$ E_H(t) = \int_{- \infty}^{\infty} w(\xi)^2 |H(x,t)|^2 d\xi $$
The function $H(.,t)$ is said to have finite weigthed energy if $E_H(t) < \infty$.
\end{def1} 

\noindent
\textbf{Remark :} A necessary condition for $E_H(0)$ to be finite is that $ \int_{-\infty}^{\infty} h_0(x) dx =0$.

\begin{prop5} 
The Cauchy problem for perturbation of Burger's viscous traveling waves has a solution $h(x,t)$ which has a holomorphic continuation belonging to $B_{\varepsilon,t}$ for all $\varepsilon \in ]0,1]$ and $t>0$ provided that $E_h(0)$ and $E_H(0)$, the weighted energy of the perturbation and of its potential are sufficiently small 
\end{prop5}

\begin{proof}
The following integral equation has to be complexified :
\begin{equation}
h(x,t) = \partial_x F_0(x,t) - \int_0^t d\tau \int_{-\infty}^{+ \infty} G(x,t-\tau;\xi) R(\tilde{f}(\xi),h(\xi,\tau)) d\xi
\end{equation}
Assuming that for each $\tau>0$, $h(x,\tau)$ has a holomorphic continuation in the strip $\Delta_{m(y_0(1-\varepsilon),\tau)}$. The integration contour in the variable $\xi$ can be shifted to the contour $\xi+ i m(y,\tau)$. The investigated $h(z,t)$ solves the resulting integral equation :
\begin{equation}
h=\partial_z F_0 + B(h)
\end{equation}
where $B$ has already been introduced in Proposition 2. 
An estimate of the linear term $\partial_z F_0$ is now required in order to get a control of its norm by some smallness condition on the initial data (here the initial weighted energies $E_h(0)$ and $E_H(0)$). 
\begin{equation} \nonumber
\begin{split}
& \partial_z F_0(z,t) =  \int_{-\infty}^{\infty} \frac{w(\xi)}{w(z)} \left( \partial_z K(z,t,\xi) - \frac{1}{2} \Phi'(f(z))\tilde{f}(z) K(z,t,\xi) \right) H_0(\xi) d\xi \\
& =  \int_{-\infty}^{\infty} \frac{w(\xi)}{w(z)} \left(  -\partial_{\xi} K(\xi,t,z)  - p'(z) K(z,t,\xi)- \frac{1}{2}\Phi'(\tilde{f}(z)) K(z,t,\xi) \right) H_0(\xi) d \xi \\
& =\frac{1}{w(z)} \int_{-\infty}^{\infty} \left( \frac{1}{2} \Phi'(\tilde{f}(\xi)) K(\xi,t,z) w(\xi) H_0(\xi) + K(z,t,\xi) w(\xi) h_0(\xi) - \right. \\
& \left. \left( \frac{1}{2} \Phi'(\tilde{f}(z)) - p'(z) \right) K(z,t,\xi) H_0(\xi) \right) d\xi
\end{split}
\end{equation}
The wave $\tilde{f(\xi)}$ being bounded wrt  $\xi \in \mathbb{R}$, using Young inequality together with the inequality :
\begin{equation} \nonumber
\| K(z,t,\xi) \|_{L^1(d\xi)} \leq \text{const} \, e^{M^2/4}
\end{equation}
we get the following estimate of $\partial_z F_0$ :
\begin{equation}
\| \partial_z F_0\|_{\varepsilon,T} \leq e^{M^2/4} \left( C_3 E_H(0) + C_4 E_h(0) + \frac{C_5 E_H(0)}{|w(iy_0(1-\varepsilon))|^{\frac{2(n-1)}{n}}}  \right)  \frac{1}{|w(iy_0(1-\varepsilon))|^{1/2} } = \alpha
\end{equation} 
Proposition 1 provides a solution if
\begin{equation} \label{init}
C \alpha p(\alpha) < 1
\end{equation}
\end{proof}
\noindent
\textit{Proof of theorem 1 :}\\
One can use the following a priori estimates for the weighted energy and conclude that initial data may actually be considered as always sufficiently small
\begin{prop6}
If $E_H(0)$  is finite, then $E_H(t)$ and $E_h(t)$ are finite for all $t>0$ and are bounded from above by $C_1 exp \left( - C_2 t \right) $.
\end{prop6}

This statement is not well-known, but it is a variation of the result of Iljin-Oleinik [\ref{IO}], its proof can be found with all details in [\ref{Xin}]. 
Now, solutions of our main theorem may be constructed through the procedure : "`wait"' until (\ref{init}) is fullfilled  denote this date $T^*$ , and apply Proposition 3 with $h(x,T^*)$ as initial data.   

\section*{References}
\begin{enumerate}
	\renewcommand\labelenumi{[\theenumi]}
	%\item \label{BB} Bardos, C.; Benachour, S. \emph{Domaine d'analycité des solutions de l'équation d'Euler dans un ouvert de $R\sp{n}$} (French)  Ann. Scuola Norm. Sup. Pisa Cl. Sci. (4)  4  (1977), no. 4, 647--687
	\item \label{Ba} Bateman, H \emph{Some recent researchs on the motion of fluids}. Monthly Weather Rev. \textbf{43} (1915), 163-170.
	\item \label{BB} Briot and Bouquet, J. Ec. Polyt. (1), cah. 36 (1856), p. 133
	\item \label{BF} Bessis, D; Fournier, JD \emph{Pole condensation and the Riemann surface associated with a shock in Burgers' equation}, J.Phys.Lett. 45 (1984), L833-L841
		\item \label{B} Burgers, JM \emph{The Nonlinear Diffusion Equation}, Dordrecht : Reidel
			%\item \label{C} Cole, Julian D. \emph{On a quasi-linear parabolic equation occurring in aerodynamics.}  Quart. Appl. Math.  9,  (1951). 225--236
		%\item \label{F} Frisch, U. \emph{Turbulence : the legacy of A.N. Kolmogorov}, Cambridge University Press, 1995
		\item \label{FT} Foias, C.; Temam, R. \emph{Gevrey class regularity for the solutions of the Navier-Stokes equations. } J. Funct. Anal.  87  (1989),  no. 2, 359--369.
		%\item \label{F} Forsyth A. R., \emph{Theory of Differential Equations} (Cambridge : Cambridge University Press),1906. 
		\item \label{Ince} Ince E.L. ; \emph{Ordinary Differential Equations}, Dover Publications, New York, 1956.
		\item  \label{IO} Iljin, A. M.; Ole\u\i nik, O. A. Asymptotic behavior of solutions of the Cauchy problem for some quasi-linear equations for large values of the time. (Russian)  Mat. Sb. (N.S.)  51 (93)  1960 191--216.

%	\item \label{doc3} Henkin, G. M.; Shananin, A. A.; Tumanov, A. E. Estimates for solutions of Burgers type equations and some applications.  J. Math. Pures Appl. (9)  84  (2005),  no. 6, 717--752
%	\item \label{doc4} Henkin, G. M. \emph{Asymptotic structure for solutions of the Cauchy problem for Burgers type equations}
	\item \label{H} Hopf, Eberhard \emph{The partial differential equation $u\sb t+uu\sb x=µu\sb {xx}$}, Comm. Pure
 Appl. Math. 3, (1950). 201--230,  
 \item \label{Gel} Gelfand, I. M. \emph{Some problems in the theory of quasi-linear equations.} (Russian) Uspehi Mat. Nauk 14 1959 no. 2 (86), 87--158. ( American Mathematical Society Translations, 33, 1963).
 \item \label{GK} Gru\`jic, Zoran; Kukavica, Igor \emph{Space analyticity for the Navier-Stokes and related equations with initial data in $L\sp p$}, J. Funct. Anal.  152  (1998),  no. 2, 447--466
 		\item \label{K} Kato, Tosio Strong $L\sp{p}$-solutions of the Navier-Stokes equation in $  R\sp{m}$, with applications to weak solutions.  Math. Z.  187  (1984),  no. 4, 471--480.
\item \label{Lax} Lax, Peter D. \emph{Weak solutions of nonlinear hyperbolic equations and their numerical computation. } Comm. Pure Appl. Math. 7, (1954). 159--193.
 \item \label{LR} Lemari\'e-Rieusset, Pierre Gilles Une remarque sur l'analyticit\'e des solutions milds des \'equations de Navier-Stokes dans $R\sp 3$. (French) [On the analyticity of mild solutions for the Navier-Stokes equations in $R\sp 3$]  C. R. Acad. Sci. Paris S\'er. I Math.  330  (2000),  no. 3, 183--186.
 \item \label{Ole} Ole\u\i nik, O. A. \emph{Uniqueness and stability of the generalized solution of the Cauchy problem for a quasi-linear equation.} (Russian) Uspehi Mat. Nauk 14 1959 no. 2 (86), 165--170 (American Mathematical Society Translations 33, 1963).
 \item \label{Sa} Sattinger, D. H. \emph{On the stability of waves of nonlinear parabolic systems}. Advances in Math. 22 (1976), no. 3, 312--355.
 %\item \label{LZ} Liu, Tai-Ping; Zumbrun, Kevin, \emph{Nonlinear stability of an undercompressive shock for complex Burgers equation}, Comm. Math. Phys. 168 (1995), no. 1, 163--186
 \item \label{Xin} Xin, Zhouping, \emph{Theory of Viscous Conservation Laws}, \emph{Some Current Topics on nonlinear Conservation laws}, AMS/IP Studies in Advanced Mathematics, Volume 15, 2000, 141-193. 
 
\end{enumerate}

\end{document}